\documentclass[12pt]{article}%
\usepackage{amsmath}
\usepackage{amsmath}
\usepackage{amssymb}
\usepackage{dsfont}
\usepackage{cite}
\newsymbol \blackbox 1004
\newcommand{\eh}{\hfill}\newlength{\sperr}

\newenvironment{proof}{{\settowidth{\sperr}{\bf\rm
Proof}%
\par\addvspace{0.3cm}\noindent\parbox[t]{1.3\sperr}
{\bf\rm P\eh r\eh o\eh o\eh f\eh }%
}}{\nopagebreak\mbox{}
$\blackbox$\par\addvspace{0.3cm}}

\def\nn{\nonumber}

\def\g{\gamma}

\def\Lam{\Lambda}
\def\s{\sigma}
\def\la{\lambda}

\def\vt{\vartheta}

\def\wt{\widetilde}

\def\p{\partial}

\def\BC{{\mathbb C}}

\def\BR{{\mathbb R}}
\def\BN{{\mathbb N}}

\def\cla{{\mathcal A}}

\def\cli{{\mathcal I}}

\def\cly{{\mathcal Y}}

\newcommand{\E}{\mathrm{e}}
\newcommand{\I}{\mathrm{i}}

\newtheorem{Pa}{Paper}[section]
\newtheorem{Tm}[Pa]{{\bf Theorem}}
\newtheorem{La}[Pa]{{\bf Lemma}}

\newtheorem{Rk}[Pa]{{\bf Remark}}

\newtheorem{Ee}[Pa]{{\bf Example}}

\title{Explicit solutions of Schr\"odinger and KdV equations
in terms of square roots of the generalised matrix eigenvalues}

\author{Alexander Sakhnovich}

\date{}
\begin{document}
\maketitle

\begin{abstract}  
In this paper, we consider matrix Schr\"odinger equation, dynamical Schr\"odinger equation and matrix KdV. 
We  construct their explicit solutions using our GBDT version of B\"acklund--Darboux transformation and square roots of the generalised matrix eigenvalues.
A separate section is dedicated to several examples including the case of strongly singular potentials.
\end{abstract}

\vspace{0.2em}


\vspace{0.3em}


\section{Introduction} \label{intro}
\setcounter{equation}{0}
Schr\"odinger and KdV equations belong to the group of the most well-known and actively
studied equations and their explicit solutions are of great interest.
In particular, B\"acklund--Darboux transformations and related dressing procedures
and commutation methods are fruitful approaches to the construction of explicit
solutions of linear and integrable nonlinear equations (see, e.g., \cite{Ci, Ge, GeT, Gu, Mi, Mar, MS, SaA2, ALS17, SaSaR, Su, Sushu, ZM} and references
therein). GBDT (generalized B\"acklund-Darboux transformation), which we use here, was first introduced in our paper \cite{SaA2}, and 
a more general version of GBDT for first order systems rationally depending on the spectral parameter
was treated in \cite{SaA3, SaSaR} (see also some references therein).

We construct GBDT and explicit solutions for the matrix Schr\"odinger equation
\begin{align}\label{I1}
-y^{\prime \prime}(x, \la)+u(x)y(x,\la)=\la y(x,\la)\quad (u=u^*), \quad y^{\prime}:=\frac{d}{d x}y,
\end{align}
for the dynamical Schr\"odinger equation
\begin{align}\label{I2}&
\I \frac{\p\psi }{\p t} (x,t)=-\frac{\p^2\psi}{\p x^2}(x,t)+u(x) \psi(x,t), \end{align}
and for the matrix KdV equation
\begin{equation}\label{I3}
\frac{ \partial  u}{ \partial t}-3 u \frac{ \partial  u}{
\partial x}-
3 \frac{ \partial  u}{ \partial x} u+ \frac{
\partial^{3}  u}{ \partial x^{3}}=0.
\end{equation}
Here, $\I$ is the  imaginary unit ($\I^2=-1$), $\la$ is the so called spectral parameter ($\la\in \BC$, $\BC$ stands for the complex plane), and $u(x)$ 
and $u(x,t)$ are  $h \times h$ matrix functions.
\begin{Rk} \label{RkD1} In \eqref{I1} and \eqref{I2}, we assume that $x$ belongs  to some finite or infinite interval $\cli$ $(x\in \cli)$,
whereas $t$ in \eqref{I2} belongs to the real axis $\BR$. Without loss of generality, we assume also that $0\in \cli$ and speak later about parameter
matrices $S(0)$ and $\Pi(0)$.
\end{Rk}
The main step in the construction of the explicit solutions of \eqref{I1}--\eqref{I3} via GBDT is the construction of the generalised
eigenfunctions $\Pi(x)$ (or $\Pi(x,t)$). We consecutively construct $\Pi$ for our three systems using square roots of the generalised matrix eigenvalues.
In this way, the results of the papers \cite{FKRS18, GKS3, KoSaTe} are further developed and wider classes of $\Pi$ and solutions of \eqref{I1}--\eqref{I3}
are obtained. 

The next section is dedicated to the general construction of the solutions of \eqref{I1}--\eqref{I3}.
Interesting examples, including the case of strongly singular potentials, are treated in Section \ref{Ex}.

As usual,  $\BN$ is the set of positive integer numbers and
$I_h$ is the $h\times h$ identity matrix.

\section{Matrix Schr\"odinger and KdV equations} \label{SK}
\setcounter{equation}{0}
{\bf 1.} GBDT for Schr\"odinger equation \eqref{I1} is determined by 3 parameter matrices.
More precisely, we choose some {\it initial} system \eqref{I1} (or, equivalently, the initial potential $u=u^*$ of Schr\"odinger equation \eqref{I1}) and fix
$n \in \BN$. Then, we fix 
$n\times n$ matrices $A$ and $S(0)$, and an $n\times m$ 
$(m=2h)$ matrix $\Pi(0)$ such that
the following relations hold:
\begin{equation} \label{SK1}
AS(0)-S(0)A^*= \Pi(0)j \Pi(0)^{*}, \quad S(0)=S(0)^*, \quad j:=\begin{bmatrix} 0 & I_h \\ - I_h &0
\end{bmatrix}.
\end{equation}
Here, $ j^*=j^{-1}=-j$. GBDT for  Schr\"odinger equation is summed up, for instance, in \cite[Sections 2,3]{FKRS18}. 
In order to construct the potentials and solutions explicitly, we set here (similar to \cite{GKS3})
$u(x)\equiv 0$. That is, our initial system is trivial.

The matrix functions $\Pi(x)$ and $S(x)$ with fixed values $\Pi(0)$ and $S(0)$ are determined by the
relations (see \cite[(3.7)]{FKRS18}):
\begin{align} & \label{SK1+}
\Pi^{\prime}(x)= A \Pi(x) \left[  \begin{array}{lr} 0 & 0 \\ I_{h} &  0 \end{array}
\right]-\Pi(x) \left[  \begin{array}{lr} 0 &
I_{h} \\  u(x) &  0 \end{array} \right],
\\ & \label{SK2-}
S^{\prime}(x)=  \Pi(x) \left[  \begin{array}{lr} 0 & 0 \\ 0 & I_{h}  \end{array}
\right]
\Pi (x)^{*} .
\end{align}
Setting $u(x)\equiv 0$ and partitioning $\Pi$ into two $n\times h$ blocks 
\begin{align}&\label{SK7-}
\Pi(x)=\begin{bmatrix}\Lam_1(x) & \Lam_2 (x)\end{bmatrix}, \quad \Pi(0)=\begin{bmatrix}\vt_1 & \vt_2 \end{bmatrix},
\end{align}
 we rewrite
the relations \eqref{SK1+} and \eqref{SK2-} as
\begin{align}& \label{SK2}
\Lam_1^{\prime}(x)=A\Lam_2(x), \quad \Lam_2^{\prime}(x)=-\Lam_1(x);
\\ & \label{SK3}
S(x)=S(0)+\int_0^x\Lam_2(\xi)\Lam_2(\xi)^*d\xi.
\end{align}
The potential and solution of the GBDT-transformed Schr\"odinger equation are expressed in terms of $\Pi(x)$ and $S(x)$ \cite{FKRS18, GKS3},
and it remains to calculate the matrix functions $\Pi(x)$ and $S(x)$ determined by \eqref{SK2}, \eqref{SK3} and the given triple $\{A,S(0), \Pi(0)\}$.
\begin{La}\label{MLa} Let  an $n\times n$ matrix $Q$ be a square root of $A$: $\, Q^2=A$. Then, the matrix functions
\begin{align}& \label{SK4}
\Lam_1(x):=-\I Q \big(\E^{\I x Q}f_1-\E^{-\I x Q}f_2\big), \quad \Lam_2(x):=\E^{\I x Q}f_1+\E^{-\I x Q}f_2,
\end{align}
where $f_1$ and $f_2$ are $n\times h$ matrices, satisfy \eqref{SK2}. Correspondingly, the matrix function  
$\Pi(x):=\begin{bmatrix}\Lam_1(x) & \Lam_2 (x)\end{bmatrix}$ satisfies \eqref{SK1+} $($where $u\equiv 0)$. In the case
\begin{align}& \label{SK4+}
-\I Q(f_1-f_2)=\vt_1, \quad f_1+f_2=\vt_2,
\end{align}
the matrix $\begin{bmatrix}\Lam_1(x) & \Lam_2 (x)\end{bmatrix}$  takes the required value $\begin{bmatrix}\vt_1 & \vt_2 \end{bmatrix}$ $($determined
by the triple  $\{A,S(0), \Pi(0)=\begin{bmatrix}\vt_1 & \vt_2 \end{bmatrix}\})$ at $x=0$.
Moreover, the integral part in \eqref{SK3} may  $($for each $Q$ and $\Pi(0))$  be explicitly calculated.
\end{La}
\begin{proof}.   Simple direct calculations show that $\Lam_1$ and $\Lam_2$ given by \eqref{SK4}
satisfy \eqref{SK2} or, equivalently, that  
\begin{align}&\label{SK7}
\Pi(x):=\begin{bmatrix}\Lam_1(x) & \Lam_2 (x)\end{bmatrix}
\end{align}
satisfies \eqref{SK1+}. Clearly, $\Pi(x)$ given by \eqref{SK4}--\eqref{SK7} takes the required value at $x=0$.
Finally, we note that the entries of $\Lam_2(x)$ are sums of the terms of the form $p_k(x)\E^{\I x c_k}$, where $c_k\in \BC$ and $p_k(x)$ are
polynomials. Hence, the last statement in the lemma is valid. 
\end{proof}
\begin{Rk}\label{RkAinv} If $\det A\not=0$ square roots $Q$ of $A$ always exist $($see, e.g., \cite[Chapter VIII, \S 6]{Gant}
with further details and references in \cite[Section 2]{ALS21}$)$. Clearly, $Q$ is invertible in this case.
Therefore, the matrices
\begin{align}
& \label{SK5}
f_1:=\left(\vt_2+\I Q^{-1}\vt_1\right)\big/2, \quad  f_2:=\left(\vt_2-\I Q^{-1}\vt_1\right)\big/2
\end{align}
are well-defined. It is immediate that $f_1$ and $f_2$ given by \eqref{SK5}  satisfy \eqref{SK4+}.  
\end{Rk}
If the solutions $Z_k$ of the matrix equations
\begin{align}
& \label{Me}
\I(QZ_1-Z_1Q^*)=f_1f_1^*, \quad \I(QZ_2+Z_2Q^*)=f_1f_2^*, \quad  -\I(QZ_3-Z_3Q^*)=f_2f_2^*
\end{align}
exist, formula \eqref{SK3} and the second equality in \eqref{SK4} yield the following representation of $S(x)$:
\begin{align}
& \label{Me2}
S(x)=\E^{\I x Q}Z_1\E^{-\I x Q^*}+\E^{\I x Q}Z_2\E^{\I x Q^*}+\E^{-\I x Q}Z_2^*\E^{-\I x Q^*}+\E^{-\I x Q}Z_3\E^{\I x Q^*}.
\end{align}

In view of \eqref{SK2}, the matrix functions $\Lam_k(x)$ also admit  an essentially less convenient than \eqref{SK4} representation
\begin{align}
& \label{Me3}
\begin{bmatrix}\Lam_1(x) \\ \Lam_2 (x)\end{bmatrix}=\E^{x\cla}\begin{bmatrix}\vt_1 & \vt_2 \end{bmatrix}, \quad \cla= \begin{bmatrix}0 & A \\
-I_n & 0 \end{bmatrix}
\end{align}
(see \cite[(3.21)]{FKRS18}).
In the case $S(0)=I_n$, explicit (although somewhat inconvenient) expressions for $S(x)$ are presented in \cite{GKS3}
in terms of the matrix exponent $\E^{\I x A_{\gamma}}$, where $A_{\gamma}$ is a $4n\times 4n$ matrix.

Assume that $y(x,\la)$ satisfies the trivial Schr\"odinger equation \eqref{I1} (with $u\equiv 0$) and put $Y_0(x,\la):=\begin{bmatrix}y(x,\la) \\ y^{\prime}( x,\la)\end{bmatrix}\in \BC^{m}$
$(m=2h)$.
Using   \cite[Proposition~3.5]{FKRS18} (for the case $u\equiv 0$) and Lemma \ref{MLa} above, we obtain the next theorem.
\begin{Tm}\label{TmSEq} Let a triple $\{A,S(0), \Pi(0)\}$ satisfy  \eqref{SK1} and  assume that $Q^2=A$.  Let the matrix functions $\Pi(x)=\begin{bmatrix}\Lam_1(x) & \Lam_2 (x)\end{bmatrix}$ and $S(x)$ be
explicitly defined by formulas \eqref{SK3} and \eqref{SK4}--\eqref{SK7}. Define $($in the points of invertibility of $S(x))$ the GBDT-transformed $h\times h$ potential 
$\wt u(x)$ by the relations
\begin{equation} \label{SK8}
\widetilde{u}(x)=2\big(X_{12}(x)+X_{21}(x)+X_{22}(x)^2\big), \quad X_{ik}(x):=\Lam_i(x)^*S(x)^{-1}\Lam_k(x).
\end{equation}
Then, the function
\begin{align}& \label{SK9}
\widetilde{y}(x, \lambda )= [I_{h} \hspace{1em}0] w_A(x, \lambda )Y_0(x, \lambda ), 
\\ & \label{SK10} w_{A}(x, \lambda ):=I_{m}- j\Pi(x)^{*}S(x)^{-1}(A- \lambda
I_{n})^{-1} \Pi(x).
\end{align}
satisfies the transformed matrix
Schr\"odinger equation
\begin{equation} \label{SK11}
-\widetilde{y}^{\, \prime \prime }(x, \lambda )+ \widetilde{u}(x) \widetilde{
y}(x, \lambda )= \lambda \widetilde{ y}(x, \lambda ).
\end{equation}
\end{Tm}
\begin{Rk} It is easy to see that the vector functions $Y_0$ have the form
\begin{equation}\label{SK18}
Y_0(x,\la)=W_0(x,\la)f_0, \quad W_0(x,\la):=\begin{bmatrix}\E^{\I x \sqrt{\la}}I_h &  \E^{-\I x \sqrt{\la}}I_h\\ \I  \sqrt{\la}\,\E^{\I x \sqrt{\la}}I_h &  -\I  \sqrt{\la} \,\E^{-\I x \sqrt{\la}}I_h\end{bmatrix}, 
\end{equation}
where $f_0\in \BC^m$ are arbitrary constant vectors. 

It is also immediate from \eqref{SK8} that $\wt u=\wt u^*$.
\end{Rk}
{\bf 2.}  The same $\Pi(x)$ and $S(x)$ provide explicit solutions of the dynamical Schr\"odinger systems
\begin{align}\label{SK19}&
\I \frac{\p}{\p t}\wt \psi(x,t)=\big(\wt H \wt \psi\big)(x,t), \quad \wt H:=-\frac{\p^2}{\p x^2}+\wt u(x).
\end{align}
Using again Lemma \ref{MLa}, we reformulate \cite[Theorem 3.1]{FKRS18} (for the case $u\equiv 0$).
\begin{Tm}\label{TmDS} Let a triple $\{A,S(0), \Pi(0)\}$ satisfy  \eqref{SK1} and assume that $Q^2=A$.  Let the matrix functions $\Pi(x)=\begin{bmatrix}\Lam_1(x) & \Lam_2 (x)\end{bmatrix}$ and $S(x)$ be
explicitly defined by formulas \eqref{SK3} and \eqref{SK4}--\eqref{SK7}. Define $($in the points of invertibility of $S(x))$ the GBDT-transformed $h\times h$ potential 
$\wt u(x)$ by the relations \eqref{SK8}.

Then, in the points of
invertibility of $S(x)$, the $m\times n$ matrix function 
\begin{align} & \label{SK20}
\wt \psi(x,t)=\begin{bmatrix} 0 & I_h \end{bmatrix}\Pi(x)^*S(x)^{-1}\E^{-\I t A}
\end{align}
satisfies the transformed  dynamical Schr\"odinger system \eqref{SK19}.
\end{Tm}
{\bf 3.}  In order to  construct explicit solutions of the matrix KdV
\begin{equation}\label{KdV}
\frac{ \partial \wt u}{ \partial t}-3\wt u \frac{ \partial \wt u}{
\partial x}-
3 \frac{ \partial \wt u}{ \partial x}\wt u+ \frac{
\partial^{3} \wt u}{ \partial x^{3}}=0,
\end{equation}
we add the variable $t$ in our matrix functions and consider 
$$\Pi(x,t)=\begin{bmatrix}\Lam_1(x,t) & \Lam_2 (x,t)\end{bmatrix}  \quad {\mathrm{and}} \quad
S(x,t) \quad (x\in\cli_1, \,\, t\in\cli_2),$$ 
where $\cli_1$ and $\cli_2$ are intervals containing $0$.
Instead of  the matrix identity \eqref{SK1}, we require
\begin{equation} \label{SK1'}
AS(0,0)-S(0,0)A^*= \Pi(0,0)j \Pi(0,0)^{*} \quad \big(S(0,0)=S(0,0)^*\big).
\end{equation}
We partition $\Pi$ into the $n\times h$ blocks 
\begin{align}&\label{SK7'}
\Pi(x,t)=\begin{bmatrix}\Lam_1(x,t) & \Lam_2 (x,t)\end{bmatrix}, \quad \Pi(0,0)=\begin{bmatrix}\vt_1 & \vt_2 \end{bmatrix}.
\end{align}
Equations \eqref{SK2} take the form
\begin{align}& \label{SK2'}
\frac{\p}{\p x}\Lam_1(x,t)=A\Lam_2(x,t), \quad \frac{\p}{\p x} \Lam_2(x,t)=-\Lam_1(x,t),
\end{align}
and another pair of PDEs is added (see \cite[p. 372]{GKS3}):
\begin{align}& \label{SK22}
\frac{\p}{\p t}\Lam_1(x,t)=4A^2\Lam_2(x,t), \quad \frac{\p}{\p t} \Lam_2(x,t)=-4A\Lam_1(x,t).
\end{align} 
Finally, $S(x,t)$ is determined by the relations (see \cite[(5.6) and (5.9)]{GKS3}):
\begin{align}& \label{SK23}
\frac{\p}{\p x}S=\Lam_2\Lam_2^*, \quad \frac{ \partial S}{ \partial t}= 4( A
\Lambda_{2}
\Lambda_{2}^{*}+
\Lambda_{2} \Lambda_{2}^{*} A^{*}+
\Lambda_{1}
\Lambda_{1}^{*}).
\end{align} 
Similar to Lemma \ref{MLa}, we derive the following lemma.
\begin{La}\label{MLa'} Let  an $n\times n$ matrix $Q$ be a square root of $A$: $\, Q^2=A$. Then, the matrix functions
\begin{align}
& \label{SK24}
 \Lam_1(x,t):=-\I Q\big(\E^{\I (x Q+4tQ^3)}f_1-\E^{-\I (x Q+4tQ^3)}f_2\big),
\\& \label{SK25}
\Lam_2(x,t):=\E^{\I (x Q+4tQ^3)}f_1+\E^{-\I (x Q+4tQ^3)}f_2,
\end{align}
where $f_1$ and $f_2$ are $n\times h$ matrices, satisfy \eqref{SK2'} and \eqref{SK22}. If \eqref{SK4+} holds,
we have  $\begin{bmatrix}\Lam_1(0,0) & \Lam_2 (0,0)\end{bmatrix}=\Pi(0,0)$. 
Moreover,  $S(x,t)$ may   be explicitly calculated $($for each $Q$ and $\Pi(0,0))$
using \eqref{SK23}.
\end{La}
\begin{Rk} \label{RkAinv'} If $\det A\not=0$, the square root $Q$ always exists 
and the matrices $f_1, f_2$ satisfying \eqref{SK4+} are given by \eqref{SK5} $($similar to the case of Remark \ref{RkAinv}$)$.
\end{Rk}
Equalities \eqref{SK2'} and the first equality in \eqref{SK23} yield
\begin{align}
& \label{SK26}
\frac{\p}{\p x}(AS-SA^*)=\frac{\p}{\p x}(\Lam_1\Lam_2^*-\Lam_2\Lam_1^*).
\end{align}
Equalities \eqref{SK22} and the second equality in \eqref{SK23} yield
\begin{align}
& \nn
\frac{\p}{\p t}(AS-SA^*)=4\big(A^2\Lam_2\Lam_2^*-\Lam_2\Lam_2^*(A^*)^2+A\Lam_1\Lam_1^*-\Lam_1\Lam_1^*A^*\big),
\\ & \nn
\frac{\p}{\p t}(\Lam_1\Lam_2^*-\Lam_2\Lam_1^*)=4\big(A^2\Lam_2\Lam_2^*-\Lam_2\Lam_2^*(A^*)^2+A\Lam_1\Lam_1^*-\Lam_1\Lam_1^*A^*\big),
\end{align}
that is,
\begin{align}
& \label{SK27}
\frac{\p}{\p t}(AS-SA^*)=\frac{\p}{\p t}(\Lam_1\Lam_2^*-\Lam_2\Lam_1^*).
\end{align}
From \eqref{SK1'}, \eqref{SK26} and \eqref{SK27} we derive
\begin{equation} \label{SK28}
AS(x,t)-S(x,t)A^*= \Pi(x,t)j \Pi(x,t)^{*}.
\end{equation}
According to the proof of \cite[Theorem 0.5]{GKS3}, relations \eqref{SK2'}--\eqref{SK23} and \eqref{SK28}
imply that  the matrix function
\begin{equation} \label{SK29}
\widetilde{u}(x,t)=2\big(X_{12}(x,t)+X_{21}(x,t)+X_{22}(x,t)^2\big),
\end{equation}
where
\begin{equation} \label{SK30}
X_{ik}(x,t):=\Lam_i(x,t)^*S(x,t)^{-1}\Lam_k(x,t),
\end{equation}
satisfies KdV \eqref{KdV}. Using also Lemma \ref{MLa'}, we obtain the following theorem.
\begin{Tm}\label{TmKdV} Let a triple $\{A,S(0,0), \Pi(0,0)\}$ satisfy  \eqref{SK1'} and  assume that $Q^2=A$.  Let the matrix functions $\Lam_1(x,t)$, $\Lam_2 (x,t)$ and $S(x,t)$ be
explicitly defined by formulas \eqref{SK23}--\eqref{SK25} and \eqref{SK4+}. Define $($in the points of invertibility of $S(x,t))$ the GBDT-transformed $h\times h$ potential 
$\wt u(x,t)$ by the relations \eqref{SK29} and \eqref{SK30}. Then, $\wt u$ satisfies KdV equation \eqref{KdV}.
\end{Tm}
\section{Examples} \label{Ex}
\setcounter{equation}{0}
Let us consider several useful examples of the potentials $\wt u$ of Schr\"odinger equations \eqref{SK11} and \eqref{SK19}
generated $($via relations \eqref{SK3}, \eqref{SK4} and \eqref{SK8}$)$ by some special triples 
$\{A,S(0), \Pi(0)=\begin{bmatrix}\vt_1 & \vt_2 \end{bmatrix}\}$ satisfying \eqref{SK1}. The corresponding
explicit solutions $\wt y$ and $\wt \psi$ follow (in terms of $\Lam_1(x)$, $\Lam_2(x)$ and $S(x)$) from Theorems \ref{TmSEq} and \ref{TmDS}, respectively.

\begin{Rk} Rational potentials $($rational extensions$)$ are of interest in applications $($see, e.g., \cite{Gran, Gran0} and references
therein$)$.  If $Q$ $($or, equivalently, $A)$ is nilpotent, it follows from \eqref{SK3}, \eqref{SK4} and \eqref{SK8} that the entries
of the potential $\wt u(x)$ are rational functions.
\end{Rk}
The simplest example is the case $A=0$.
\begin{Ee}\label{Ee2} \cite[p. 371]{GKS3}. Assume that 
$A=Q=0$ and $\vartheta_{1}=0$. Then, \eqref{SK1} holds for any $S(0)=S(0)^*$ and any $\vt_2$.
The equalities \eqref{SK4+} are valid in the case $f_1+f_2=\vt_2$. Thus, by virtue of  \eqref{SK3}, \eqref{SK4} and \eqref{SK8} we have
\begin{align}&\nn
\Lambda_{1}(x)=0, \hspace{1em}
\Lambda_{2}(x)= \vartheta_{2}, \hspace{1em}
S(x)= S(0)+x \vartheta_{2}  \vartheta_{2}^{*},
\\ &\label{3.1}
\wt u(x)=2\big( \vartheta_{2}^{*}\big(
S(0)+x \vartheta_{2}  \vartheta_{2}^{*}\big)^{-1}
\vartheta_{2}\big)^{2}
\end{align}
\end{Ee}
\begin{Rk} The expression \eqref{3.1} for $\wt u$ may be simplified $($especially for the scalar case $h=1)$  in an easy way, see \eqref{3.1'} below.
\end{Rk}
Indeed, assume that $S(0)$ is invertible and rewrite \eqref{3.1} as
 \begin{align} &\label{Z1}
\wt u(x)=2\big( \vartheta_{2}^{*}\big(
I_n+x \theta  \vartheta_{2}^{*}\big)^{-1}
\theta\big)^{2}, \quad \theta:=S(0)^{-1}\vt_2.
\end{align}
Then, using geometric progressions, we rewrite  the resolvent $\big(I_n+x \theta  \vartheta_{2}^{*}\big)^{-1}$ in the form
 \begin{align} \nn
\big(
I_n+x \theta  \vartheta_{2}^{*}\big)^{-1}
&=
I_n-x\theta\left(\sum_{k=1}^{\infty}(-x\vt_2^*\theta)^{k-1}\right)\vt_2^*
\\ &\label{Z2}
=I_n-x\theta(I_h+x\vt_2^*\theta)^{-1}\vt_2^*.
\end{align}
The equality \eqref{Z2} holds for small $x$ and so (in view of the analyticity) for all points of invertibility. Taking into account  \eqref{Z2}, we obtain
 \begin{align} \nn
\vartheta_{2}^{*}\big(
I_n+x \theta  \vartheta_{2}^{*}\big)^{-1}
\theta&=\vt_2^*\theta+(I_h+x\vt_2^*\theta)^{-1}\vt_2^*\theta-(I_h+x\vt_2^*\theta)(I_h+x\vt_2^*\theta)^{-1}\vt_2^*\theta
\\ & \label{Z3}
=(I_h+x\vt_2^*\theta)^{-1}\vt_2^*\theta.
\end{align}
Relations \eqref{Z1} and \eqref{Z3} yield
 \begin{align}  &\label{3.1'}
\wt u(x)=2\big((I_h+x\vt_2^*\theta)^{-1}\vt_2^*\theta\big)^2.
\end{align}

Our next example deals with a slightly more complicated subcase of the case $A=0$.
\begin{Ee}\label{Ee3}  Let $h=1$, $n=2$, $A=0$ and
\begin{align}
& \label{3.2}
Q=\begin{bmatrix} 0 & 1 \\ 0 & 0\end{bmatrix}, \quad \vt_1=\begin{bmatrix} b \\ 0\end{bmatrix}, \quad \vt_2=\begin{bmatrix} c \\ 0\end{bmatrix}, \quad S(0)=\begin{bmatrix} 0 & 0 \\ 0 & d\end{bmatrix} \quad (d\not=0), 
\end{align}
where $b,c,d \in \BR$. Clearly, we have $Q^2=A=0$. 
\end{Ee}
In the case of the given above $\vt_1$ and $\vt_2$, we also have $\Pi(0)j\Pi(0)^*=0$. Thus, \eqref{SK1} holds for any $S(0)=S(0)^*$ 
and we choose $S(0)$ as in \eqref{3.2}. For $f_i=\begin{bmatrix} f_{i1} \\ f_{i2}\end{bmatrix}$ $(i=1,2)$ relations \eqref{SK4+} are equivalent to
\begin{align}
& \label{3.2+}
f_{11}+f_{21}=c, \quad f_{12}=-f_{22}=\I b/2.
\end{align} 
Further we assume that \eqref{3.2+} (and so \eqref{SK4+}) is valid and use Theorems \ref{TmSEq} and \ref{TmDS}.
Since $Q^2=0$, the series representations of $\E^{\pm \I  x Q}$ and relations \eqref{SK4} and \eqref{SK4+} imply that
\begin{align}
& \label{3.3}
\Lam_1(x)=-\I Q(f_1-f_2)=\vt_1, \quad \Lam_2(x)=f_1+f_2+\I Qx(f_1-f_2)=\vt_2-x\vt_1.
\end{align}
Using \eqref{SK3} and taking into account \eqref{3.2} and \eqref{3.3}, we obtain
\begin{align}
& \label{3.4}
S(x)=\begin{bmatrix} \g(x) & 0 \\ 0 & d\end{bmatrix}, \quad \g(x):=(b^2/3)x^3-bcx^2+c^2x.
\end{align}
Finally, relations \eqref{SK8}, \eqref{3.3} and \eqref{3.4} yield
\begin{align}
\nn
\wt u(x)&=\frac{4b}{\g(x)}(c-b x)+\frac{2}{\g(x)^2}(c-b x)^4
\\ & \label{3.5}
=\frac{2(bx-c)}{\g(x)^2}\big((b^3/3)x^3-b^2cx^2+bc^2x-c^3\big).
\end{align}
\begin{Rk}\label{StSing}  The case of strong singularities of $\wt u(x)$ at $x=0$, for instance,
\begin{align}
& \label{3.6}
\wt u(x) \sim \ell(\ell+1)/x^2 \quad {\mathrm{for}} \quad x\to 0,
\end{align}
is of special interest \cite{Gran0, KoSaTe0, KoSaTe}.
Formula \eqref{3.5} $($for the Example \ref{Ee3} above$)$ shows that  
\begin{align}
& \label{3.7}
\wt u(x)=\frac{2}{x^2}\big(1+O(x)\big) \quad {\mathrm{for}} \quad c\not=0, \,\, x\to 0; \\
& \label{3.7'}
 \wt u(x)=\frac{6}{x^2}  \quad {\mathrm{for}} \quad c=0, \quad b\not=0.
\end{align}
Thus, we have the case $\ell=1$ for $c\not=0$ and $\ell=2$ for $c=0$. 
\end{Rk}
 Another example for the case $\ell=1$ and an example for the case $\ell =3$
have been treated in \cite{KoSaTe0} and \cite{KoSaTe}, respectively. (In both cases, we had  $S(0)=0$ but $\det A \not=0$.)
\begin{Ee}\label{Ee3.6}
Let us simplify formulas for fundamental solutions  in  our example \eqref{3.7'} where $\ell=2$.
Since $h=1$ and $c=0$, relations \eqref{3.2}--\eqref{3.4} and \eqref{SK10} yield $A=0$,  $\vt_2=0$ and
\begin{align}
& \label{3.8}
\Pi(x)=\begin{bmatrix} b & -bx \\ 0 & 0\end{bmatrix}, \quad S(x)=\begin{bmatrix}  (b^2/3)x^3 & 0 \\ 0 & d\end{bmatrix}, 
\\ & \label{3.9}
\begin{bmatrix} 1 & 0 \end{bmatrix}w_A(x,\la)=\begin{bmatrix} 1-\frac{3}{\la x^2} & \frac{3}{\la x} \end{bmatrix}.
\end{align}
Hence, Theorem \ref{TmSEq} and formulas \eqref{SK9} and \eqref{SK18} imply that
\begin{align}
& \label{3.10}
\phi(x,\la):=\begin{bmatrix} 1 & 0 \end{bmatrix}w_A(x,\la)W_0(x,\la)\begin{bmatrix} 1 \\ 0 \end{bmatrix}=\E^{\I x\sqrt{\la}}\left(1+\frac{3\I}{\sqrt{\la}\, x}-\frac{3}{\la x^2}\right)
\end{align}
satisfies Schr\"odinger equation with the potential $\wt u(x)=6/x^2:$
\begin{align}
& \label{3.10+}
-\wt y^{\,\prime \prime}(x, \la)+\frac{6}{x^2}\wt y(x,\la)=\la \wt y(x,\la).
\end{align}
 In view of \eqref{SK9} and \eqref{SK18}, another solution $\chi(x,\la)$ of this
equation is obtained by the substitution of $-\sqrt{\la}$ instead of $\sqrt{\la}$ on the right-hand side of \eqref{3.10}, that is,
\begin{align}
& \label{3.11}
\chi(x,\la)=\E^{-\I x\sqrt{\la}}\left(1-\frac{3\I}{\sqrt{\la}\, x}-\frac{3}{\la x^2}\right).
\end{align}
Clearly, the branch of $\sqrt{\la}$ in \eqref{3.10} and \eqref{3.11} may be chosen in an arbitrary way
$($the same branch for both formulas$)$. 

Now, it is easy to construct a nonsingular at $x=0$  solution $\cly$ of  the Schr\"odinger equation \eqref{3.10+}
as a linear combination of $\phi$ and $\chi$. Namely, we take $\cly(x,\la):=\phi(x,\la)-\chi(x,\la)$. In view of \eqref{3.10}
and \eqref{3.11}, it is easily verified that
\begin{align}
\nn
\cly(x,\la)=&(1+\I x\sqrt{\la})\left(1+\frac{3\I}{\sqrt{\la}\, x}-\frac{3}{\la x^2}\right)
\\ &  \nn
-(1-\I x\sqrt{\la})\left(1-\frac{3\I}{\sqrt{\la}\, x}-\frac{3}{\la x^2}\right)+O(1) 
=\frac{6\I}{\sqrt{\la}\, x}-\frac{6\I}{\sqrt{\la}\, x}+O(1)
\end{align}
for $x\to 0$. That is, $\cly(x,\la)$ is nonsingular at $x=0$.
\end{Ee}
Using double commutation method, S. Albeverio, R. Hryniv, and Ya. Mykytyuk \cite{AHM} studied
the change of $\ell$ in the term $(\ell/x)\s_1$ in a  radial Dirac system when an
eigenvalue is removed or inserted. (See also a related paper \cite{T}.)
The change of $\ell$ in the case of GBDT for radial Dirac systems was studied in \cite[pp. 237--239]{SaSaR}
(see also the references therein) and the analog of this result for the Schr\"odinger equations would be of
interest.

{\bf Acknowledgments}  {This research    was supported by the
Austrian Science Fund (FWF) under Grant  No. Y-963.}

\end{document}